\documentclass{amsart} 

\usepackage{amsmath,amsthm}     
\usepackage{amssymb}            
\usepackage{euscript}           
\usepackage{enumerate,calc}
\usepackage[matrix,arrow,curve,frame]{xy}    

\xymatrixcolsep{1.9pc}                          
\xymatrixrowsep{1.9pc}
\newdir{ >}{{}*!/-5pt/\dir{>}}                  

\newtheorem{thm}[subsection]{Theorem}
\newtheorem{defn}[subsection]{Definition}
\newtheorem{prop}[subsection]{Proposition}

\newtheorem{cor}[subsection]{Corollary}
\newtheorem{lemma}[subsection]{Lemma}

\theoremstyle{definition}  
\newtheorem{example}[subsection]{Example}

\newtheorem{remark}[subsection]{Remark}

  {\end{list}}

\newcommand{\dfn}{\textbf} 

\newcommand{\mdfn}[1]{\dfn{\mathversion{bold}#1}} 

\newcommand{\tens}              {\otimes}               
\newcommand{\iso}               {\cong}  

\newcommand{\cat}{\EuScript}    
\newcommand{\cC}{{\cat C}}

\newcommand{\cM}{{\cat M}}
\newcommand{\cN}{{\cat N}}

\newcommand{\cP}{{\cat P}}

\newcommand{\cT}{{\cat T}}
\newcommand{\cU}{{\cat U}}

\newcommand{\Top}{{\cat Top}}
\newcommand{\Gtop}{\Top(G)}
\newcommand{\Zttop}{\Top(\Z_2)}

\newcommand{\sSet}{s{\cat Set}}

\newcommand{\field}[1]  {\mathbb #1} 
\newcommand{\A}         {\field A}
\newcommand{\R}         {\field R}

\newcommand{\Z}         {\field Z}
\newcommand{\C}         {\field C}

\DeclareMathOperator*{\colim}{colim}
\DeclareMathOperator*{\hocolim}{hocolim}

\DeclareMathOperator{\spec}{Spec}
\DeclareMathOperator{\Hom}{Hom}

\DeclareMathOperator{\Map}{Map}

\DeclareMathOperator{\sk}{sk}
\DeclareMathOperator{\Sk}{Sk}

\DeclareMathOperator{\Ex}{Ex}

\DeclareMathOperator{\sd}{sd}
\DeclareMathOperator{\cosk}{cosk}

\newcommand{\ra}{\rightarrow}                   
\newcommand{\lra}{\longrightarrow}              
\newcommand{\llra}[1]{\stackrel{#1}{\lra}}      

\newcommand{\mapiso}{\llra{\cong}}
\newcommand{\we}{\llra{\sim}}                   

\newcommand{\cof}{\rightarrowtail}              

\newcommand{\trcof}{\stackrel{\sim}{\cof}}

\newcommand{\inc}{\hookrightarrow}              

\newcommand{\blank}{-}                          
\newcommand{\id}{id}                            

\newcommand{\norm}[1]{|#1|}       

\newcommand{\non}{\text{non}}

\newcommand{\Cech}{\check{C}}
\newcommand{\CCech}{\v{C}ech\ }

\newcommand{\scop}{{\Delta^{op}}}

\newcommand{\sing}{Sing}
\newcommand{\Sing}{\sing\,}

\newcommand{\assign}{\mapsto}

\newcommand{\ovcat}{\downarrow}

\newcommand{\bd}[1]{\partial\Delta^{#1}}

\newcommand{\adjoint}{\rightleftarrows}

\newcommand{\del}[1]{\Delta^{#1}}

\newcommand{\he}{\simeq}

\newcommand{\Sm}{Sm}

\newcommand{\rea}[1]{|{#1}|}             
\newcommand{\map}{\rightarrow}

\newcommand{\ceck}[1]{\Cech(#1)}         
\newcommand{\oceck}[1]{\Cech^{o}(#1)}    
\newcommand{\Spc}{\text{Spc}}            
\newcommand{\fcup}{\cup}
\newcommand{\Gal}{\text{Gal}}
\newcommand{\stn}{s_{\leq n}}

\numberwithin{equation}{section}

\begin{document}

\title{Hypercovers in Topology}

\author{Daniel Dugger}
\author{Daniel C. Isaksen}

\address{Department of Mathematics\\ Purdue University\\ West
Lafayette, IN 47907 } 

\address{Department of Mathematics \\ University of Notre Dame\\ 
Notre Dame, IN 46556}

\email{ddugger@math.purdue.edu}
\email{isaksen.1@nd.edu}

\thanks{
The second author was supported by an NSF Postdoctoral Research Fellowship.}

\date{November 25, 2001}

\begin{abstract}
We show that if $U_*$ is a hypercover of a topological
space $X$ then the natural map $\hocolim U_* \ra X$ is a weak
equivalence.  This fact is used to construct topological realization
functors for the $\A^1$-homotopy theory of schemes over real and
complex fields.
\end{abstract}

\maketitle


\section{Introduction}

Let $X$ be a topological space, and let $\cU = \{U_a\}$ be an open
cover of $X$.  From this data one may build the \mdfn{\CCech complex
$\Cech(\cU)_*$}, which is the simplicial space
\[ \xymatrix{
\coprod U_{a_0}
& \coprod U_{a_{0}a_{1}} \ar@<0.5ex>[l]\ar@<-0.5ex>[l]
& \coprod U_{a_{0}a_{1}a_{2}} \ar@<0.6ex>[l]\ar[l]\ar@<-0.6ex>[l]
\cdots
}
\]
Here $U_{a_{0} \cdots a_{n}} = U_{a_{0}} \cap \cdots \cap U_{a_{n}}$,
and the face maps are obtained by omitting indices---we have chosen
not to draw the degeneracies for typographical reasons.  Segal
\cite{S1} proved that if $X$ has a partition of unity subordinate to
$\cU$ then the map $\rea{\Cech(\cU)_*} \ra X$ is a {\it homotopy\/}
equivalence, where $\rea{\blank}$ denotes geometric realization.  Our
first goal in this paper is to generalize this result to the following
theorem.

\begin{thm}
\label{th:cech}
For every open cover $\cU$ of $X$, the natural map $\hocolim
\Cech(\cU)_* \ra X$ is a weak equivalence.
\end{thm}

There are two steps in the argument.  First, we prove that
$\rea{\Cech(\cU)_*} \ra X$ is a weak equivalence for arbitrary open
covers.  It is possible to deduce this from Segal's result, making use
of the fact that weak equivalences are detected by spheres, and
spheres always have partitions of unity.  But instead of going this
route we give a proof that avoids Segal's theorem completely, and is
quite elementary.

The second step is to deal with the difference between
$\rea{\Cech(\cU)_*}$ and $\hocolim \Cech(\cU)_*$.  For any simplicial
object $W_*$ in a model category, there are general criteria for when
its geometric realization agrees with its homotopy colimit
(cf. \cite[Th.~19.6.4]{H}); unfortunately these criteria apply only
when the objects $W_n$ are all cofibrant, and we are definitely not
assuming that the open sets $U_a$ and their intersections are
cofibrant.  To get around this we prove a curious theorem (given in
Appendix A) that when computing homotopy colimits for topological
spaces one never has to worry about this cofibrancy issue.  Strange,
but true.

The main goal of this paper is generalizing Theorem~\ref{th:cech} so
that it applies to `hypercovers', rather than just \CCech covers.
These are defined in detail in Section \ref{se:hypercover}, but for
now we will just give an intuitive definition.  An open hypercover of a
space $X$ is a simplicial space $U_*$ such that
\begin{enumerate}[(1)]
\item Each $U_n$ is a disjoint union of open subsets of $X$,
\item The spaces appearing in $U_0$ are an open cover of $X$,
\item The spaces in $U_1$  cover the double
intersections of those in level 0,
\item The spaces in $U_2$ cover the triple intersections of those in
level 1, and so on.
\end{enumerate}
Of course making sense of (4)---especially the `and so on'
part---requires a certain amount of bookkeeping, which is why we are
postponing the formal definition.  But the essence is that hypercovers
are like \CCech complexes except that instead of taking the double
intersections at level 1 we may refine them further, and we may
continue this refining process at each level.  Our second main result
is then

\begin{thm}
\label{th:hypercover}
If $U_*$ is an open hypercover of a space $X$, then the natural map
$\hocolim U_* \ra X$ is a weak equivalence.
\end{thm}

This result could almost be considered folklore since everyone
immediately agrees it's true, but a proof seems to be missing from the
literature.  One might consider tackling it by appealing to the
Whitehead theorem, proving an isomorphism on fundamental groupoids and
homology with local coefficients.  This is the approach taken in
\cite[Prop.~8.1]{F} in the related context of \'etale hypercovers, but
this is messy and obscures in computation the underlying geometric
explanation of the theorem.  In the case of topological spaces, the
isomorphism on fundamental groupoids was the subject of the paper
\cite{RT} (although they only dealt with \CCech complexes, not
hypercovers).  The approach we take here, on the other hand, is very
elementary. The idea is to reduce to the case of \CCech covers in a
clever way.

Our interest in these results arose from attempts to understand
topological realization functors in the $\A^1$-homotopy theory of
schemes \cite{MV}.  Given an algebraic variety $X$ defined over $\C$,
there is an associated topological space $X(\C)$ obtained by giving
$X$ the analytic topology.  Of course this should extend to a map of
`homotopy theories' from the Morel-Voevodsky category $\Spc(\C)$ to
the category of topological spaces.  In \cite{MV} this extension is
only provided at the level of homotopy categories, but we are
interested in extending it to the model category level.  The key fact
needed to make this work is precisely Theorem~\ref{th:hypercover}.
This is worked out in detail in Section~\ref{se:MV}, following the
basic program of \cite{I} (also outlined in \cite[Rem.~8.2]{D2}).
We also prove that taking analytic spaces for schemes
defined over $\R$ induces a Quillen map from
$\Spc(\R)$ to $\Z_2$-equivariant topological spaces.

Finally, we give in this paper several interesting corollaries to
Theorem~\ref{th:hypercover}.  On the whole these seem too disparate to
recount in the introduction, but as an example let us 
mention two of them.  We refer the reader to
Sections~\ref{se:hocolims} and \ref{se:hypercover} for more results
like these.

\begin{cor}
\label{co:covspace}
Let $E\ra B$ be any map which is locally split (for example, a covering space),
and form the associated \CCech complex $\Cech(E)_*$ given by
\[ \Cech(E)_n := E^{n+1}_B = E\times_B E \times_B \cdots \times_B E. 
\]
Then the natural map $\hocolim \Cech(E)_* \ra B$ is a weak equivalence.
\end{cor}

\begin{cor}
Let $\cU$ be an open cover of a space $X$ with the property that
every finite intersection $U_{a_0\cdots a_n}$ is covered by other
elements of $\cU$.  Form the diagram consisting of all the $U_a$'s
and all the inclusions between them.  Then the homotopy colimit of
this diagram is weakly equivalent to $X$.
\end{cor}

The first corollary is an immediate consequence of Proposition
\ref{pr:genhyp}, and the second is restated and proved as Proposition
\ref{pr:compcover}(c).

Using open covers to give homotopy decompositions for spaces, or to
detect weak equivalences, is of course a classical topic.  In addition
to \cite{S1} it is worthwhile to mention \cite{Mc1}, \cite{Mc2}, and
\cite{Dk}.  Hypercovers were invented by Verdier in \cite[Expose V,
Sec.~7]{SGA4}, where they were used as a way of computing sheaf
cohomology in arbitrary Grothendieck topologies.

\medskip

We would like to express our thanks to Bill Dwyer, Phil Hirschhorn,
Michael Mandell, and Jeff Smith 
for several useful conversations about these results.

\subsection{Notation, terminology, and other annoyances}

We assume that the reader is familiar with homotopy colimits, and in a
few places also with the theory of model categories.  The original
reference for the latter is \cite{Q}, but we generally follow \cite{H}
in notation and terminology (\cite{Ho} is also a good reference).
Regarding homotopy colimits, \cite{H} uses `$\hocolim D$' to denote
the result of applying a certain explicit formula to any diagram $D$.
This has the disadvantage that the resulting object has the correct
homotopy type only when the diagram consists entirely of cofibrant
objects.  We instead adopt the position that `$\hocolim D$' should
{\it always\/} denote the correct homotopy-invariant construction: it
is obtained by first applying cofibrant-replacement to the objects in
the diagram, and only then using the usual explicit formulas.  In
model-theoretic terms, homotopy colimit is the left derived functor of
the ordinary colimit functor, when the category of diagrams is given
the projective model structure (see below).

Having made the previous point, we now get to say that for topological
spaces it isn't really necessary.  This is definitely a non-standard
fact, but we've banished it to Appendix A so it won't distract the
reader from the general theme of the paper. On the other hand, it is a
useful result and we'd like to call the reader's attention to it: when
taking homotopy colimits for diagrams of topological spaces one
doesn't first have to make all the spaces involved cofibrant.  The
usual formulas are already homotopy-invariant.

\label{se:proj}
We review one last piece of machinery, used often in the body
of the paper.  Given a small category $I$, recall that there is a
model structure on the category of diagrams $\sSet^I$ such that a map
is a weak equivalence (resp., fibration) if it is so in every spot of
the diagram \cite[Sec.~13.8]{H}.  We call this the {\it projective\/}
model structure on $\sSet^I$, and the cofibrant diagrams have the
property that the homotopy colimit and ordinary colimit are weakly
equivalent.

Finally, some notation: Throughout this paper our open covers
$\cU=\{U_a\}$ are always indexed by a set $A$.  In particular, we are
allowing the possibility that $U_a=U_{a'}$ for different values $a\neq
a'$.  For every finite set $\sigma = \{ a_0, \ldots, a_n \}$ in $A$,
we'll write $U_\sigma$ or $U_{a_0 \cdots a_n}$ for $U_{a_0} \cap
\cdots \cap U_{a_n}$.  Also, once and for all we fix our model for
$\del{n}$ as the subset of $\R^{n+1}$ consisting of $(n+1)$-tuples $t
= (t_{0}, \ldots, t_{n})$ such that $0 \leq t_{i} \leq 1$ for all $i$
and $\Sigma_{i = 0}^{n} t_{i} = 1$.  The symbol $\Top$ denotes the
category of {\it all\/} topological spaces---we don't assume any
hypotheses like compactly-generated.


\section{\CCech complexes}
\label{se:cech}

The purpose of this section is to prove the following:

\begin{thm}
\label{th:segalwe}
For any open cover $\cU$ of a topological space $X$, the natural map
$\pi: \rea{\Cech(\cU)_*} \ra X$ is a weak equivalence.
\end{thm}

We start by recalling the following result and its corollary:

\begin{prop}[Gray]
\label{pr:Gray}
Let $f\colon X\ra Y$ be a map of spaces and let $U$ and $V$ form an
open cover of $Y$.  Suppose that the induced maps
\[ f^{-1}U \ra U, \qquad f^{-1}V \ra V, \quad \text{and} \quad
 f^{-1}(U\cap V)
\ra U\cap V \]
are all weak equivalences.  Then $X\ra Y$ is also a weak equivalence.
\end{prop}

This is proven (in more generality) in \cite[16.24]{Gr}, using an
elegant small-simplices argument.  With enough technology it can also
be done by a Whitehead-type theorem: it's easy to see that $X\ra Y$ is
an isomorphism on $\pi_0$, a souped-up van Kampen theorem yields the
isomorphism on $\pi_1$, and for homology with local coefficients one
uses the Mayer-Vietoris exact sequence.  Gray's argument is much
nicer, though.

\begin{cor}[May]
\label{co:May}
Let $f\colon X \ra Y$ be a map of spaces and let $\cU = \{U_a\}$ be an open
cover of $Y$.  
Suppose that $f^{-1}U_{\sigma} \ra U_{\sigma}$ is a
weak equivalence for every finite set $\sigma$ of indices.  Then $X\ra
Y$ is also a weak equivalence.
\end{cor}

May deduces the generalization by a quick application of Zorn's Lemma
\cite[Cor.~1.4]{M2}: look at the set of all opens $W$ such that
$f^{-1}(W\cap U_\sigma) \ra W\cap U_\sigma$ is a weak equivalence for
all $\sigma$, including $\sigma=\emptyset$.  This set has a maximal
element, and Gray's result shows it must be $X$.  In an earlier paper
McCord proved a more general version of this result \cite[Th.~6]{Mc1},
but the proof is quite a bit more complicated.

\begin{proof}[Proof of Theorem~\ref{th:segalwe}]
Given any open set $V$ in $X$, the space $\pi^{-1}(V)$ is homeomorphic to
the space $\rea{\Cech(\cU')_*}$, where $\cU'$ is the open cover
$\{U_a \cap V\}$ of the space $V$.  This definitely uses the fact that
$V$ is open.

We want to consider the maps $\pi^{-1}(U_\sigma) \ra U_\sigma$, but in
this case the cover $\cU'$ of $U_\sigma$ actually contains the whole
space $U_\sigma$ as one of its elements.  From the following lemma we
know that under this condition $\rea{\Cech(\cU')_*} \ra U_\sigma$ is a
weak equivalence; so by Corollary \ref{co:May} the map
$\rea{\Cech(\cU)_*} \ra X$ is a weak equivalence as well.
\end{proof}

\begin{lemma}
\label{le:splitcover}
Let $\cU$ be an open cover of $X$ such that $U_b=X$ for
some index $b$.  Then the natural map $\rea{\Cech(\cU)_*} \ra X$ is a weak
equivalence (in fact, a homotopy equivalence). 
\end{lemma}

\begin{proof}
There is a section $\chi\colon X \ra \rea{\Cech(\cU)_*}$ obtained from
the map $U_b\tens \del{0} \ra \rea{\Cech(\cU)_*}$ and the identification
$U_b=X$.  We only need to show that $\chi\pi$ is homotopic to the
identity.

Let $\Cech(\cU)_*\times I$ be the simplicial space obtained by crossing
all the levels of $\Cech(\cU)_*$ with the unit interval.  Then
$\rea{\Cech(\cU)_* \times I}$ is the quotient
\[ \Biggl [\coprod_{a_0\cdots a_n} U_{a_0\cdots a_n} \times \del{n} \times I
\Biggr ] / \sim
\]
where the relations are the usual ones, not affecting the $I$ factor
at all.  Define a map $\rea{\Cech(\cU)_* \times I} \ra \rea{\Cech(\cU)_*}$ in
the following way.  Take an element
$(x,t_0,\ldots,t_n,s)$ where $x$ belongs to $U_{a_0\cdots
a_n}$ and $(t_0,\ldots,t_n)$ belongs to $\del{n}$, and send it to the element
$(x,1-s,st_0,\ldots,st_n)$ in the factor $U_{ba_0\ldots a_n} \tens
\del{n+1}$.  
This definition respects the various identifications.

Now, there is also an obvious map $f\colon \rea{\Cech(\cU)_* \times I}
\ra \rea{\Cech(\cU)_*} \times I$ induced by sending $(x,t,s)$ to $((x,t),s)$.
We claim that this is a homeomorphism, thereby giving us a homotopy
$\rea{\Cech(\cU)_*}\times I \ra \rea{\Cech(\cU)_*}$ between $\chi\pi$ and the
identity.  The reason $f$ is a homeomorphism is just because geometric
realization and crossing with $I$ are both left adjoints, and the
right adjoints are easily seen to commute.  It is important that $I$
and $\Delta^n$ are locally compact Hausdorff so that the relevant
mapping spaces with compact-open topologies have the correct
adjointness properties.
\end{proof}

\subsection{Connection with Segal's results}
\label{se:segalret}

To close this section we make the connection between our
Theorem~\ref{th:segalwe} and the result proven in \cite{S1}.  Segal
doesn't explicitly deal with \CCech complexes, but the objects he
deals with turn out to be homeomorphic to them.  This connection will
be needed later on.

Let $A$ be the indexing set for a cover $\cU$.  We have already
introduced the \CCech complex $\Cech(\cU)_*$, but if $A$ is given an
ordering we may also consider the \dfn{ordered \CCech complex}
$\oceck{\cU}_*$ which is often easier to work with.  This is the
simplicial space given by $\oceck{\cU}_{n} = \coprod_{a_{0} \cdots
a_{n}} U_{a_{0} \cdots a_{n}}$, where the coproduct ranges over all
{\em ordered} multi-indices in $A$.  That is, we only consider
multi-indices for which $a_{0} \leq a_{1} \leq \cdots \leq a_{n}$.
Note that there is an inclusion of simplicial spaces $\oceck{\cU}_* \ra
\ceck{\cU}_*$.

\begin{prop}
\label{pr:ceck=oceck}
The map
$\oceck{\cU}_* \ra \ceck{\cU}_*$
induces a homotopy equivalence
$\rea{\oceck{\cU}_*} \ra \rea{\Cech(\cU)_*}$. 
\end{prop}

\begin{proof}
For any (not necessarily ordered) multi-index $a_{0} \cdots a_{n}$,
there is a canonical reordering $a_{\sigma 0} \cdots a_{\sigma n}$
such that $a_{\sigma 0} \leq \cdots \leq a_{\sigma n}$.  If
$a_{i} = a_{j}$ for some $i < j$, then always choose $\sigma i < \sigma j$.
This allows us to define an inverse map $\rea{\Cech(\cU)_*} \map
\rea{\oceck{\cU}_*}$.
  If $(x, t)$ is an element of $U_{a_{0} \cdots a_{n}} \otimes \del{n}$,  
then send $(x,t)$ to the element $(x, \sigma
t)$ of $U_{a_{\sigma 0} \cdots a_{\sigma n}} \otimes \del{n}$,
where $\sigma t$ is defined by $(\sigma t)_{i} = t_{\sigma i}$.

One composition is the equal to the identity.  It remains to construct
a homotopy $H: \rea{\Cech(\cU)_*} \times I \map \rea{\Cech(\cU)_*}$ between the
other composition and the identity.  As in the proof of Lemma
\ref{le:splitcover}, we use the space $\rea{\ceck{\cU}_* \times I}$
rather than $\rea{\Cech(\cU)_*} \times I$.  We define $H$ as follows: An
element $(x, t)$ of $U_{a_{0} \cdots a_{n}} \otimes \del{n}$ is
equivalent in $\rea{\Cech(\cU)_*}$ to the element $(x, t_{0}, \ldots, t_{n},
0, \ldots, 0)$ of $U_{a_{0} \cdots a_{n} a_{\sigma 0} \cdots a_{\sigma
n}} \otimes \del{2n+1}$.  Also, $(x, \sigma t)$ is equivalent in
$\rea{\Cech(\cU)_*}$ to the element $(x, 0, \ldots, 0, t_{\sigma 0}, \ldots,
t_{\sigma n})$ of $U_{a_{0} \cdots a_{n} a_{\sigma 0} \cdots a_{\sigma
n} } \otimes \del{2n+1}$.  Define $H((x,t), s)$ to be the element
\[
(x, st_{0}, \ldots, st_{n}, (1-s)t_{\sigma 0}, \ldots, (1-s)t_{\sigma n})
\]
of
$U_{a_{0} \cdots a_{n} a_{\sigma 0} \cdots a_{\sigma n}} 
  \otimes \del{2n+1}$.
\end{proof}

\begin{prop}
\label{pr:oceck=sdceck}
Let $\cU$ be an open cover of a space $X$ indexed by a set $A$.
Consider
the realization of the simplicial space
\[ [n] \assign \coprod_{\sigma_0 \subseteq \cdots \subseteq \sigma_n}
U_{\sigma_n},
\] 
where the coproduct is indexed by chains of nonempty, finite subsets
of $A$.  This realization is homeomorphic to the realization
$\rea{\oceck{\cU}_*}$ of the ordered \CCech complex and is homotopy equivalent
to $\rea{\Cech(\cU)_*}$.
\end{prop}

The realization in the above proposition 
is the object considered in \cite{S1}.  
The ordered \CCech complex is another construction of the same space, which
for us
seems somewhat easier to work with.  One disadvantage, of course, is
that it is not natural: a total ordering on $A$ must be chosen to
begin with.

\begin{proof}
The second claim follows from the first claim and 
Proposition \ref{pr:ceck=oceck}.

For the first claim,
it is convenient to use a slightly unusual construction of
$\rea{\oceck{\cU}_*}$.  When forming the
geometric realization, instead of forming Cartesian products with
$\del{k}$ we instead form products with $\sd \del{k}$; since they are
homeomorphic it doesn't matter which one we use.  Given this, the key
observation is that we can coordinatize $\sd \del{k}$ in the following
way: assuming that the vertices of $\del{k}$ are labelled by the
numbers $0,\ldots,k$ in the usual way, a point on $\sd \del{k}$ is
represented uniquely by a chain of proper inclusions $\sigma_0\subset
\cdots \subset \sigma_j$ of subsets of $\{0,\ldots,k\}$ together with
an element $t$ of $\del{j}$.  Essentially, the chain of subsets
determines in which sub-simplex the point lies, and then $t$ gives
local coordinates inside that sub-simplex.

Using this coordinate scheme, we can write down maps
in both directions between the
two realizations
\[ \Biggr [\coprod_{\sigma_0 \subseteq \cdots \subseteq \sigma_n}
U_{\sigma_n} \times \del{n}
\Biggl ]/\sim 
\qquad
\text{and} \qquad
\Biggr [\coprod_{a_0\leq \cdots \leq a_k}
U_{a_0\cdots a_k} \times \sd\del{k}
\Biggl ]/\sim. 
\]
For instance, let's give the map from left to right.  
Using degeneracy relations, 
a point $p$ in
the left space can be represented by a chain of {\it proper\/}
inclusions $\sigma_0 \subset \cdots \subset \sigma_n$, a point $x$ of
$U_{\sigma_n}$, and an element $t$ of $\del{n}$.  Let
$a_0, a_1,\ldots,a_k$ be the ordered list of elements of $\sigma_n$.  
The chain
$\sigma_*$ together with $t$ defines a point $s$ in $\sd
\del{k}$, and so we map $p$ to the pair $(x,s)$.  It is easy to see
that this map is well-defined and continuous, and just as easy to
write down its inverse.
\end{proof}

In the case that $\{U_a\}$ admits a partition of unity $\{\psi_a\}$ it
is fairly easy to see that the map $\pi\colon \rea{\oceck{\cU}_*} \ra X$ admits
a section: First, a point $x$ of $X$ has a neighborhood which intersects
the support of $\psi_a$ only for finitely many indices
$a=a_0,\ldots,a_n$.  The section $\chi$ sends $x$ to the point of
$\rea{\oceck{\cU}_*}$ represented by $(x,t)$ in $U_{a_0\cdots a_n} \tens \del{n}$
where $t_i=\psi_{a_i}(x)$.  One has to check that $\chi$ is continuous
(use the local-finiteness of the partition of unity), and that
$\chi\pi \he \id$ via a straight-line homotopy.  See Proposition 4.1
of \cite{S1}.


\section{Passing to homotopy colimits}
\label{se:hocolims}

The results of the previous section all concerned geometric
realizations.  In this section we translate these into results about
various homotopy colimits.  In general, there is a `Reedy cofibrancy'
condition on simplicial spaces which guarantees that geometric
realization and homotopy colimit agree.  Unfortunately our \CCech
complexes are not Reedy cofibrant, due to the fact that the
open sets appearing in them are not necessarily cofibrant spaces.
However, Theorem \ref{th:hocolim} shows that in the category of
topological spaces this cofibrancy issue is unimportant: homotopy
colimits can be computed naively, without first making things
cofibrant.  This fact saves the day.

\bigskip

\begin{thm}
\label{th:seghoco}
If $\cU$ is an open cover of a space $X$, then the natural map
$\hocolim \Cech(\cU)_* \ra X$ is a weak equivalence.
\end{thm}

\begin{proof}
By Theorem \ref{th:hocolim}, we can compute the homotopy colimit in
the Strom model category.  In this model structure the \CCech complex
is Reedy cofibrant (it has free degeneracies in the sense of
Definition \ref{de:freedeg}), and so the realization already has the correct
homotopy type.  Theorem~\ref{th:segalwe} now gives the result.
\end{proof}

Here are several alternative formulations:

\begin{prop}
\label{pr:hocostuff}
Let $A$ be an indexing set for the cover $\cU$, and let $\cP_A$ denote
the partially ordered set consisting of all nonempty finite subsets of
$A$.  Let $\Gamma$ denote the functor $\cP_A^{op} \ra \Top$
which sends $\sigma$ to $U_\sigma$.  Then the natural map $\hocolim
\Gamma \ra X$ is a weak equivalence.
\end{prop}

\begin{proof}
To construct $\hocolim \Gamma$ we can take the realization
of the simplicial replacement for $\Gamma$ (by Theorem \ref{th:hocolim}
we don't need to first make the spaces cofibrant).  
That is, we take
the realization of the simplicial space
\[ [n] \assign \coprod_{\sigma_0 \subseteq \cdots \subseteq \sigma_n}
U_{\sigma_n},
\] 
where the coproduct is indexed by chains of nonempty, finite subsets
of $A$.  Now Proposition \ref{pr:oceck=sdceck} tells us that this
realization is homotopy equivalent to $\rea{\Cech(\cU)_*}$, so Theorem
\ref{th:segalwe} finishes the proof.
\end{proof}

\begin{cor}
\label{co:hocostuff1}
Let $\cP_{\cU}$ denote the subcategory of $\Top$ whose objects are the
open sets $U_a$ belonging to $\cU$ together with their finite
intersections; the morphisms are the inclusions of open subsets of
$X$.  Let $\Gamma$ denote the inclusion functor $\cP_{\cU} \ra \Top$.
Then the natural map $\hocolim \Gamma \ra X$ is a weak equivalence.
\end{cor}

\begin{proof}
Consider the obvious functor $F\colon \cP_A^{op} \ra \cP_{\cU}$
sending $\sigma$ to $U_\sigma$.  We will show that it is homotopy
cofinal, so pick an object $V$ in $\cP_{\cU}$ and look at the
undercategory $(V\ovcat F)$.  It suffices to show that any map $K\ra
N(V\ovcat F)$ can be extended over the cone on $K$, as $K$ ranges over
all finite simplicial sets.  Every $n$-simplex $s$ in $K$ maps to a
chain of open sets $V \ra U_{\sigma_0} \ra U_{\sigma_1}\ra \cdots \ra
U_{\sigma_n}$ in $(V\ovcat F)$.  Since $K$ has only finitely-many
non-degenerate simplices, only finitely-many of the $U_\sigma$ will
ever appear.  Define $\mu$ to be the union of all the $\sigma_i$
arising from the map $K\ra N(V\ovcat F)$.  To extend the map over
$CK$, we send the cone on $s$ to the $(n+1)$-simplex corresponding to
the chain $V \ra U_\mu \ra U_{\sigma_0} \ra U_{\sigma_1}\ra \cdots \ra
U_{\sigma_n}$.
\end{proof}

The following corollary was shown to us by Bill Dwyer.  Let
$(\Top\ovcat X)_{\cU}$ denote the full subcategory of $(\Top\ovcat X)$
consisting of all maps $Z\ra X$ that factor through 
the space $E = \coprod_a U_a$.
Let $\Gamma\colon (\Top\ovcat X)_{\cU} \ra \Top$ be the canonical functor
sending $Z\ra X$ to $Z$.  We would like to claim that 
the homotopy colimit of the diagram $\Gamma$ is weakly equivalent to $X$,
but $(\Top\ovcat X)_{\cU}$ is not a small category.  So we choose an
infinite cardinal $\kappa$ larger than the size of $E$
and restrict to the spaces $Z$ that have at most $\kappa$ elements.
As the proof of the corollary indicates, the weak homotopy type of 
$\hocolim \Gamma$ is independent of the choice of $\kappa$, as long as
$\kappa$ is sufficiently large.

\begin{cor}
\label{co:hocostuff2}
For the functor $\Gamma\colon (\Top\ovcat X)_\cU \ra \Top$ defined
above, the natural map $\hocolim \Gamma \ra X$ is a weak equivalence.
\end{cor}

\begin{proof}
The $n$th level
of the \CCech complex is $E_X^n:=E\times_X E \times_X \cdots \times_X
E$ ($n$ factors).  Let's write $\cC =(\Top\ovcat X)_{\cU}$, for brevity.
So we have the functor $F\colon\scop \ra \cC$ given by $[n] \assign
E^n_X$.  The composition $\scop \ra \cC \ra \Top$ is just $\Cech(\cU)_*$.
Because of Theorem \ref{th:seghoco},
it will be enough to show that $F$ is homotopy cofinal.

For this we pick an object $z:Z\ra X$ in $\cC$ and show that $(z\ovcat
F)$ is contractible.  This undercategory is isomorphic to the category
of simplices of $K$, where $K$ is the simplicial set sending $[n]$ to
$\Hom_{\cC}(z,E^n_X)$.  But observe that $\Hom_{\cC}(z,E^n_X)$ is
equal to $T^n$ where $T=\Hom_{\cC}(z,E_X)$.  So $K$ is the simplicial
set $[n] \mapsto T^n$, which is contractible because $T$ is nonempty
(using the fact that $z\colon Z\ra X$ factors through $E$).  Thus
$(z\ovcat F)$ is isomorphic to the category of simplices of a
contractible simplicial set, and therefore has a contractible nerve.
\end{proof}

\begin{cor}[Small simplices theorem]
Let $\Sing_{\cU} X$ denote the simplicial set whose $n$-simplices are
the maps $\del{n} \ra X$ that factor through some $U_a$.  Then
$\Sing_{\cU} X \ra \Sing X$ is a weak equivalence.
\end{cor}

\begin{proof}
Let $\cP_A$ be the category defined in Proposition~\ref{pr:hocostuff},
where $A$ is the indexing set for the cover.  Consider the diagram
$\Gamma\colon \cP_A^{op} \ra \sSet$ defined by
$\Gamma(\sigma)=\Sing(U_\sigma)$.  By general nonsense $\rea{\hocolim
\Gamma} \he \hocolim \rea{\Gamma}$.  Also, there is a commutative diagram
\[
\xymatrix{
\hocolim_{\cP_A^{op}} \rea{\Sing U_\sigma} \ar[r] \ar[d] &
  \rea{\Sing X} \ar[d] \\
\hocolim_{\cP_A^{op}} U_\sigma \ar[r] & X   }
\]
in which the vertical maps are weak equivalences because the natural
map $\rea{\Sing Y} \map Y$ is a weak equivalence for every space $Y$.
We know from Proposition \ref{pr:hocostuff} that the bottom horizontal
map is a weak equivalence, so the top horizontal map is also a weak
equivalence.  We conclude that the map $\hocolim \Gamma \ra \Sing X$
is a weak equivalence of simplicial sets.  Therefore, we shall
compare $\hocolim \Gamma$ and $\Sing_{\cU} X$.

For the moment, assume that $A$ is finite.  Notice that $\cP_A^{op}$
is a Reedy category \cite[Def.~5.2.1]{Ho}, where we think of all the
maps as being directed upward.  Since there are no non-identity
downward maps, the fibrations are objectwise in the Reedy model
structure on $\sSet^{\cP_A^{op}}$ (see \cite[Th.~5.2.5]{Ho}).  So in
this case the Reedy and projective model structures (cf. Section
\ref{se:proj}) are the same.  In particular, a Reedy-cofibrant diagram
is also projective-cofibrant, which guarantees that the homotopy
colimit and the ordinary colimit are weakly equivalent.

The functor $\Gamma$ may be checked to be Reedy cofibrant: at the spot
indexed by $\sigma = \{a_0,\ldots,a_n\}$, the latching object is the subobject
of $\Sing U_{\sigma}$ consisting of all simplices which are
contained in some other $U_b$.  The fact that it is actually a
subobject says that the latching map is a cofibration.  
So we know that $\hocolim \Gamma$ and $\colim \Gamma$ are weakly
equivalent.
It is easy to check that $\colim
\Gamma \iso \Sing_{\cU} X$.  We have shown that if $\cU$ is a finite
cover, then $\Sing_{\cU} X$ is weakly equivalent to $\Sing X$.

Now let $A$ be arbitrarily large.  For any finite subcollection
$\cU'$, let $\cup{\cU'}$ denote the union of the open sets in $\cU'$.
Then we know the map
$\Sing_{\cU'} (\cup{\cU'}) \ra \Sing(\cup{\cU'})$ is a weak
equivalence.  But $\Sing_{\cU} X \ra \Sing X$ is the filtered colimit
of these maps, where the indexing category is the poset of all finite
subcollections $\cU'$. This uses that each space $\Delta^n$ is 
compact.
Our result now follows from the fact that
filtered colimits of simplicial sets preserve weak equivalences.
\end{proof}


\section{Hypercovering Theorems}
\label{se:hypercover}

In this section we define hypercovers, and then prove our main result,
Theorem~\ref{th:hypercover}.  We go on to deduce various corollaries.

\medskip

Before giving a rigorous definition of hypercovers, we need to recall
a few pieces of machinery related to simplicial objects.  For any
category $\cC$, let $s\cC$ denote the category of simplicial objects
in $\cC$.  Likewise, let $\stn\cC$ denote the category of truncated
simplicial objects of dimension $n$.  There is the obvious forgetful
functor $\sk_n\colon s\cC \ra \stn\cC$, and if $\cC$ has all finite
limits then $\sk_n$ has a right adjoint called $\cosk_n$; these are
the \dfn{skeleton} and \dfn{coskeleton} functors.  If $U_*$ belongs
to $s\cC$ then we'll
often abbreviate $\cosk_n(\sk_n U)_*$ as just $\cosk_n U_*$.  Finally, the
\mdfn{$n$th matching object $M_n U$} is defined to be the $n$th object
of $\cosk_{n-1} U_*$.  There is a canonical map of simplicial spaces $U_*
\ra \cosk_{n-1} U_*$, and in level $n$ it gives $U_n \ra M_n U$.
In levels less than $n$, this map is the identity.
We write $\cosk^X_n$ for the $n$th
coskeleton functor for $s(\Top\ovcat X)$.

These definitions have somewhat easier interpretations when $\cC$ is
the category of topological spaces.  To describe these, note that
any simplicial set may be regarded as a simplicial space which is
discrete in every dimension, and if $U_*$ and $W_*$ are
simplicial spaces then the set of maps from $U_*$ to $W_*$ has a natural
topology coming from the compact-open topology on function
spaces.  Using these observations, one checks that 
\begin{enumerate}[(i)]
\item $U_n \iso \Map(\del{n},U_*)$,
\item $[\cosk_n U]_k \iso \Map( sk_n \del{k} , U_*)$, and
\item $M_n U \iso \Map (\bd{n}, U_*)$.
\end{enumerate} 
The first property is immediate from the Yoneda lemma.  The second
property follows from the first and the adjunction between
$\sk_n$ and $\cosk_n$.  The third property is a special case
of the second.

Finally, say that a map of spaces $Z \ra X$ is an \dfn{open covering
map} if it is isomorphic to a map of the form $\coprod_a U_a \ra X$
where $\{ U_a \}$ is an open cover of $X$.

\begin{defn}  A \dfn{hypercover} of a space
$X$ is an augmented simplicial space $U_* \ra X$ such that the maps
$U_n \ra M^X_n U$ are open covering maps for all $n\geq 0$.  Here
$M^X_n U$ denotes the $n$th matching object of $U_*$ computed in the
category $s(\Top\ovcat X)$ of simplicial spaces over $X$.  
\end{defn}

Note that $M^X_0 U \iso X$, so the condition for $n=0$ says that $U_0
\ra X$ is an open covering map.  Also $M^X_1 U\iso U_0 \times_X U_0$,
so when $n=1$ we are requiring $U_1 \ra U_0 \times_X U_0$ to be an
open covering map.  The reader should be aware that when $n>1$ 
the objects $M_n U$ and $M^X_n U$ turn out to be isomorphic, so one
can forget about the extra complication of the overcategory.

Using properties (i)--(iii) above, it can be checked that 
if $U_* \map X$ is a hypercover and $K \map L$ is an inclusion
of finite simplicial sets, then the map
$\Map(L, U_*) \map \Map(K, U_*)$ is also an open covering map. 
From this, it follows that
$\cosk^X_n U_* \map X$
is a hypercover whenever 
$U_* \map X$ is a hypercover.  Also, each map
$U_k \map [ \cosk^X_n U ]_k$ is an open covering map.

We leave it to the reader to check that in a hypercover each $U_n$
must be a disjoint union of open subsets of $X$, and that \CCech
complexes are the hypercovers for which the maps $U_n \ra M^X_n U$
are all {\it isomorphisms\/}.  Generalizing this, a hypercover $U_*
\ra X$ is called \dfn{bounded} if there exists an $N$ such that the
maps $U_n \ra M^X_n U$ are isomorphisms for all $n> N$.  The smallest
such $N$ for which this happens is called the dimension of the
hypercover.  Said intuitively, the bounded hypercovers of dimension
$N$ are the hypercovers for which the refinement process stops after
the $N$th level.  A hypercover $U_* \map X$ has dimension at most $N$ if and 
only if $U_* \iso \cosk^X_N U_*$.

\begin{lemma}
\label{le:fd-real}
If $U_* \ra X$ is a bounded hypercover, then $\hocolim{U_*} \ra X$ is a
weak equivalence.
\end{lemma}

A more detailed version of the following proof, given in the context
of an arbitrary Grothendieck topology, appears in \cite{DHI}.

\begin{proof}
We proceed by induction, starting from the fact that bounded
hypercovers of dimension $0$ are just \CCech covers and
therefore are handled by Theorem~\ref{th:seghoco}.

Suppose that $U_* \ra X$ is a bounded hypercover of dimension $n+1$.
Define $V_*$ to be $\cosk_n U_*$, so $V_*$ is a bounded hypercover of
dimension at most $n$.  Therefore, we may assume by induction that
$\hocolim{V_*} \map X$ is a weak equivalence.  The canonical map $U_*
\ra V_*$ gives an open covering map $U_{n+1} \ra V_{n+1}$, by the very
definition of what it means for $U_*$ to be a hypercover (since
$V_{n+1}=M_{n+1}U$).  In fact, one can check that $U_k \ra V_k$ is an
open covering map for all $k$.

Consider the following bisimplicial object, augmented horizontally by
$V_*$:

\[ \xymatrix{
 V_* & U_* \ar[l] & U_* \times_{V_*} U_* \ar@<0.5ex>[l]\ar@<-0.5ex>[l] &
 \cdots \ar@<0.6ex>[l]\ar[l]\ar@<-0.6ex>[l] }
\]
The $k$th row is the (augmented) \CCech complex for the open covering map
$U_k \ra V_k$.  Note that for $0\leq k \leq n$ the $k$th row
is the constant simplicial object with value $U_k$
because $U_k \ra V_k$ is the
identity.  Call this bisimplicial object (without the horizontal
augmentation) $W_{**}$.  

Let $D_*$ denote the diagonal of $W_{**}$.  Standard homotopy theory
tells us that $\hocolim{D_*}$ may be computed (up to weak equivalence)
by first taking the homotopy colimits of the rows of $W_{**}$, and
then taking the homotopy colimits of the resulting simplicial object.
But the homotopy colimit of the $k$th row is just $V_k$ by Theorem
\ref{th:seghoco}.  Since $V_*$ is a bounded hypercover of dimension
at most $n$, we have assumed that $\hocolim{V_*}$ is weakly equivalent to $X$.
So $\hocolim{D_{*}} \ra X$ is a weak equivalence.

We claim that $U_*$ is a retract, over $X$, of $D_*$.  Note first that
one has, in complete generality, a map $U_* \ra D_*$; in dimension $k$
it is the unique horizontal degeneracy $W_{0k} \ra W_{kk}$.

To produce a map $D_* \ra U_*$ it is enough to give $\sk_{n+1} D_* \ra
\sk_{n+1} U_*$, because $U_*=\cosk_{n+1} U_*$.  Notice that $\sk_n D_*
= \sk_n U_*$.  Choosing any face map $[0] \ra [n+1]$ gives a map
$W_{n+1,n+1} \ra W_{0,n+1}$, which is just $D_{n+1}\ra U_{n+1}$.  This
induces a corresponding map $\sk_{n+1} D_* \ra sk_{n+1} U_*$ as
desired.

It is straightforward to check that $U_* \ra D_* \ra U_*$ is the
identity (because $U_*=\cosk_{n+1} U_*$ one only has to check it on
$(n+1)$-skeleta), and all the maps commute with the augmentations down
to $X$.  
We have already shown that $\hocolim{D_*} \map X$ is a weak equivalence.
Since $\hocolim{U_*} \map X$ is a retract of $\hocolim{D_*} \map X$, 
it must also be a weak equivalence.
\end{proof}

\begin{thm}
\label{th:hypmain}
If $U_* \ra X$ is a hypercover then the maps $\hocolim U_* \ra
\rea{U_*} \ra X$ are all weak equivalences.
\end{thm}

\begin{proof}
The fact that $\hocolim U_* \ra \rea{U_*}$ is a weak equivalence
follows just as in Theorem \ref{th:seghoco} for the case of \CCech
complexes: we may compute the homotopy colimit in the Strom model
category, where the simplicial object $U_*$ is Reedy cofibrant since
it has free degeneracies (Definition~\ref{de:freedeg}).

To show that $\rea{U_*} \ra X$ is a weak equivalence, note first that
we have an isomorphism $\pi_k\rea{U_*} \ra \pi_k \rea{\cosk_{k+1} U_*}$.
 (This is true for any map of simplicial spaces $X_* \ra Y_*$
which is an isomorphism on $(k+1)$-skeleta---an easy proof is to
apply the singular functor everywhere to get into bisimplicial sets,
then use the diagonal in place of realization.)
But $\cosk_{k+1} U_*$ is a bounded hypercover, so Lemma
\ref{le:fd-real} tells us that $\pi_k \rea{\cosk_{k+1} U_*} \mapiso
\pi_k X$.
\end{proof}

\subsection{Complete covers}

In this section we don't quite consider hypercovers, but rather a
related concept which captures the same phenomena.  This second
approach was suggested to us by Jeff Smith.

\begin{defn}
An open cover $\cU = \{ U_a \}$ of a space $X$ is called \dfn{complete}
if for all finite sets $\sigma$ of indices, the intersection
$U_\sigma$ is covered by elements of $\cU$.
It is called a \dfn{\CCech cover} if every $U_\sigma$ is
again an element of the cover.
\end{defn}

Complete covers appear in \cite[Satz~2.2]{DT}, where they were used in the
context of identifying quasi-fibrations.  The paper \cite{Mc1} then used
them to detect weak equivalences.

We blur the distinction between a cover and the full subcategory that
it spans inside the category of open sets of $X$.  Given a cover
$\cU$, we can construct an associated simplicial space in the
following way: For any $n \geq 0$, let $P_{n}$ denote the category of
nonempty subsets of $\{0,\ldots,n\}$, where the maps are the
inclusions.  Note that the assignment $[n]\assign P_n$ defines a
cosimplicial category in the obvious way. (Application of the nerve
functor everywhere gives the cosimplicial space $[n] \ra \sd\del{n}$.)

Define $\Omega_*$ to be the simplicial space
\[ [n] \assign \coprod_{F\colon P_n^{op} \ra \cU} F(\{0,\ldots,n\}),
\]
where the coproduct runs over all functors $P_n^{op} \ra \cU$.  
The faces and degeneracies are induced by those in $P$ in the
expected way.

To give a point in $\Omega_3$, for example, is to give the following
data: 
\begin{enumerate}[(1)]
\item
A sequence of opens $U_0,\ldots,U_3$ in $\cU$, 
\item $6$ open
subsets
$U_{01},U_{02},\ldots,U_{23}$ in $\cU$ such that $U_{ij}\subseteq
U_i\cap U_j$;
\item $4$ open subsets $U_{012},\ldots,U_{123}$ in $\cU$ such that
$U_{ijk}\subseteq U_{ij}\cap U_{jk} \cap U_{ik}$;
\item An open subset $U_{0123}$ in $\cU$ which is contained in all the
$U_{ijk}$;
\item A point on $U_{0123}$.
\end{enumerate}
It is usually helpful to think of these open sets as indexed by the faces
of a $3$-simplex.  

In forming the \CCech complex of a cover $\cU$ we are throwing in all
the finite intersections $U_\sigma$ into the higher levels of the
simplicial object, and these are typically objects which are not in
$\cU$ itself.  The simplicial object $\Omega_*$ is in some sense the closest
thing we can get to a \CCech complex while requiring all the open sets
to belong to $\cU$.

\begin{prop}
\label{pr:compcover}
\mbox{}\par
\begin{enumerate}[(a)]
\item If the cover $\cU$ is complete then $\Omega_*$ is a hypercover
of $X$.
\item Regarding $\cU$ as a category,
let $\Gamma\colon \cU \ra \Top$ be the obvious
inclusion.  Then $\hocolim \Gamma\he\norm{\Omega_*}$.
\item If the cover $\cU$ is complete then the natural map $\hocolim
\Gamma \ra X$ is a weak equivalence.
\end{enumerate}
\end{prop}

\begin{proof}
For part (a), consider the full subcategory $\bar{P}_n$ of $P_n$
consisting of all objects except for $\{ 0, 1, \ldots, n \}$.  Then
the matching space $M_n \Omega$ is equal to
\[ \coprod_{\bar{F}\colon \bar{P}_n^{op} \ra \cU}\,\, 
   \Bigl [ \bigcap_{\sigma \in \bar{P}_n} \bar{F}(\sigma)
   \Bigr ].
\]
For example, a point in $M_3 \Omega$ 
is determined by the data in (1)--(3)
above, together with a point in $U_{012}\cap U_{013}\cap U_{023} \cap
U_{123}$.  

Since the cover is complete, for each functor $\bar{F}\colon
\bar{P}_n^{op} \ra \cU$ and each element $x$ of $\cap_{\sigma \in
\bar{P}_n} \bar{F}(\sigma)$, there exists an extension $F$ of $\bar{F}$
to $P_n^{op}$ such that $x$ belongs to $F(\{ 0, \ldots, n \})$.  This
shows that $\Omega_n \map M_n \Omega$ is an open covering map, which
finishes part (a).

Part (b) is almost trivial, given the right machinery.  To form
$\hocolim \Gamma$ we can work in the Strom model structure on $\Top$
(see Appendix A), where we first take the simplicial replacement
\[ [n] \assign \coprod_{U_0\ra \cdots \ra U_n} U_0 
\]
and then form the realization.  Here the coproduct is indexed over
all functors $\del{n} \ra \cU$, where $\del{n}$ denotes the category
of $n$ composable maps.  Note that $\Omega_*$ was formed in almost the
same way as the simplicial replacement of $\Gamma$, except we indexed
the coproduct by functors $P_n^{op} \ra \cU$.  Each $P_n$ is
essentially just a subdivision of $\del{n}$, so it's not
surprising that $\norm{\Omega_*}$ is another model of the homotopy
colimit.

In somewhat more detail: Let $\sd'$ denote the `opposite' of the usual
subdivision functor on $\sSet$, in which the orientations of all the
simplices have been changed so that they point away from the
barycentres, rather than towards them.  (We need this because we are
using $P_n^{op}$ rather than $P_n$.)  The functor $\sd'$ has a right
adjoint $\Ex'$.  There is a natural `first vertex map' $\sd' K \ra K$,
inducing $K \ra \Ex' K$.  Given our diagram $\Gamma\colon \cU \ra
\Top$, the realization of the simplicial replacement is isomorphic to
the coend $\Gamma \tens_{\cU} B$, where $B \colon \cU \ra \sSet$ sends
$U_a$ to the classifying space $B(U_a \ovcat \cU)$.  Likewise, one
checks that the realization of $\Omega_*$ is isomorphic to the coend
$\Gamma \tens_\cU \Ex' B$, where $\Ex' B$ is the obvious composite functor.
The natural map $B \ra \Ex' B$ is an objectwise weak equivalence.  
The object $B$ of $\sSet^{\cU}$ is cofibrant
(see \cite[Cor.~15.8.8]{H}),
where this diagram category has the projective model
structure described in Section~\ref{se:proj}.  
The exact same
arguments show that $\Ex' B$ is also cofibrant in this structure.  So
we have an objectwise weak equivalence between two cofibrant diagrams.
The diagram $\Gamma\colon \cU \ra \Top$ is objectwise cofibrant (since
we are working with the Strom model structure on $\Top$), and so by
\cite[Cor.~19.3.5]{H} it follows that $\Gamma\tens_\cU B \ra
\Gamma\tens_\cU \Ex' B$ is a weak equivalence.

Finally, part (c) is an immediate consequence of (a), (b), and
Theorem~\ref{th:hypmain}.
\end{proof}

The following corollary was originally proven by McCord
\cite[Th.~6]{Mc1}, but is an easy consequence of our hypercovering
theorem.  It generalizes May's result from Corollary~\ref{co:May},
which handled the case of \CCech covers.  For the proof we will need
the following observations: (1) If $U \ra X$ is an open covering map
and $f\colon Y\ra X$ is any map, there is an induced open covering map
$Y\times_X U \ra Y$.  (2) If $U_* \ra X$ is a hypercover and
$f\colon Y\ra X$ is a map of spaces, one gets a hypercover
$f^{-1} U_* \ra Y$ whose space in level $n$ is $Y\times_X U_n$.

\begin{cor}
Let $f\colon X\ra Y$ be a map of spaces.  Suppose there is a complete
cover $\cU = \{U_a\}$ of $Y$ such that each $f^{-1}(U_a) \ra U_a$ is a weak
equivalence.  Then $f$ itself is a weak equivalence.
\end{cor}

\begin{proof}
From $\cU$ form the associated hypercover $\Omega_*^Y$ as described in
the paragraph preceding Proposition~\ref{pr:compcover}.  Pulling this
back to $X$ gives a hypercover $\Omega_*^X:=f^{-1}\Omega_*^Y$, as
described above (note that this is {\it not\/} the hypercover
associated to the covering $\{f^{-1}U_a\}$).  Now $f$ induces a map
$\Omega_*^X \ra \Omega_*^Y$ compatible with the augmentations.  This
map of simplicial spaces is a levelwise weak equivalence, by
assumption.  Upon taking homotopy colimits we get
\[
 \xymatrix{
 \hocolim \Omega_*^X \ar[r]^{\sim} \ar[d]_{\sim} & \hocolim \Omega_*^Y
 \ar[d]^{\sim} \\
  X \ar[r] & Y,}
\]
and so we conclude that $X\ra Y$ is also a weak equivalence.
\end{proof}

\subsection{Generalized hypercovers for topological spaces}

Up until now we have only considered open covers, but now we turn to a
broader notion.  We'll say that a map $p\colon E\ra B$ of spaces is a
\dfn{generalized cover} if it is locally split: that is, every element
of $B$ has a neighborhood $U$ such that $p^{-1}(U) \ra U$ admits a
section.
Observe that covering spaces, and in fact fibre bundles in general,
are generalized covers.  The point for us is that generalized covers
and open covers generate the same Grothendieck topology on topological
spaces.

\begin{defn}
An augmented simplicial space $U_* \ra X$ is a \dfn{generalized
hypercover} of $X$ if the maps $U_n \ra M^X_n U$ are generalized covers.
\end{defn}

\begin{prop}
\label{pr:genhyp}
If $U_*$ is a generalized hypercover of $X$ then $\hocolim U_* \ra X$
is a weak equivalence. 
\end{prop}

\begin{proof}
Results from \cite{DHI}, in the context of an arbitrary Grothendieck
topology, show that this is a consequence of Theorem~\ref{th:hypmain}.
The essential point is that generalized hypercovers can all be refined
by open hypercovers.

To deduce this from the results of \cite{DHI} we do the following:
Pick a regular cardinal $\lambda$ larger than the size of all the sets
in $U_*$.  Let $\Top_\lambda$ denote the category of topological
spaces of size less than $\lambda$, and make it into a Grothendieck
site via the usual notion of open cover.  Form the universal model
category $U(\Top_\lambda)$ (see the paper \cite{D2}) and localize it
with respect to the set $S$ consisting of all maps $\hocolim V_* \ra
X$, where $V_*\ra X$ is an open hypercover.  Theorem~\ref{th:hypmain}
implies that that there is a `realization map' $U(\Top_\lambda)/S \ra
\Top$.  The results of \cite{DHI} say that in $U(\Top_\lambda)/S$ one
actually knows that $\hocolim U_* \ra X$ is a weak equivalence for all
{\it generalized\/} hypercovers $U_*$, and then applying our
realization functor tells us this must hold in $\Top$ as well.
\end{proof}

Corollary \ref{co:covspace} is an immediate consequence of the above
proposition.

\begin{example}
Let $G$ be a topological group and consider the usual covering space
$\xi\colon EG\ra BG$.  Form the \CCech complex $\ceck{\xi}_*$, which is
a generalized hypercover of $BG$.  Using only the fact that $EG$ has a
free $G$-action, one can see that the $n$th level of $\ceck{\xi}_*$ is
homeomorphic to $G^n \times EG$, and the face and degeneracy maps are
the familiar ones of the two-sided bar construction $B(*,G,EG)$.  Now
using that $EG$ is contractible, we find that $\ceck{\xi}_*$ is
levelwise weakly equivalent to the simplicial space
\[ \xymatrix{
{*}
& G \ar@<0.5ex>[l]\ar@<-0.5ex>[l]
& G\times G \ar@<0.6ex>[l]\ar[l]\ar@<-0.6ex>[l]
\cdots.
}
\]
The above proposition tells us that $\rea{\ceck{\xi}_*} \he BG$, and so
in this way we recover the usual bar construction for $BG$.
\end{example}


\section{Topological realization functors for $\A^1$-homotopy theory}
\label{se:MV}

Let $k$ be a field.  Morel and Voevodsky \cite{MV} produced a
model category $\Spc(k)$ which captures the `motivic homotopy theory'
of smooth schemes over $k$.  Here $\Spc(k)$ stands for `spaces over
$k$'.  It is the category of simplicial presheaves on the
Nisnevich site of smooth schemes over $\spec k$.  

When $k$ comes with an embedding $k \inc \C$, then any $k$-scheme $X$
gives rise to a topological space $X(\C)$ consisting of its
$\C$-valued points with the analytic topology.  A natural expectation
is to use this functor to relate $\Spc(k)$ to the usual model category
$\Top$ of topological spaces.  Morel and Voevodsky showed how to extend
this functor on the level of homotopy categories (by somewhat awkward
methods), but they didn't produce functors at the model category level.  In
this section we use Proposition~\ref{pr:genhyp} to produce such
functors, with the small provision that we have to replace $\Spc(k)$
with a Quillen-equivalent variant.  We also address the situation when
$k\inc \R$, in which case one can construct topological realization
functors into $\Z_2$-equivariant spaces.

As in \cite{D2}, a Quillen pair $L\colon \cM \adjoint \cN \colon R$
will be called a {\it Quillen map\/} $\cM\ra\cN$.

\medskip

\subsection{The Complex case}

Let $\cT$ denote either the Zariski, Nisnevich, or \'etale
Grothendieck topology on the category $\Sm/k$ of smooth $k$-schemes.
In the terminology of \cite{D2}, let $\Spc'(k)_{\cT}$ denote the
universal model category built from $\Sm/k$ subject to the following
relations:
\begin{enumerate}[(1)]
\item $X\amalg Y \we (X\fcup Y)$ (here $\amalg$ denotes the coproduct
in our model category, whereas $\fcup$ denotes disjoint union of schemes);
\item $\hocolim U_* \we X$ for any $\cT$-hypercover $U_*$ of a smooth
scheme $X$ (called `basal hypercovers' in \cite{DHI});
\item $X\times \A^1 \we X$.
\end{enumerate}
(Relation (1) is morally a special case of (2), but must be included
separately for technical reasons---see \cite{DHI}).

The model categories $\Spc(k)_{\cT}$ and $\Spc'(k)_{\cT}$ have the
same underlying category and the same class of weak equivalences, but
differ in their notions of cofibration and fibration.  They are
injective and projective versions of the same homotopy theory.

\begin{thm}
\label{th:Creal}
There are Quillen maps $\Spc'(k)_{et} \ra \Top$ and
$\Spc'(k)_{Nis} \ra \Top$
sending a smooth $k$-scheme
$X$ to $X(\C)$.  
\end{thm}

\begin{proof}
By general nonsense from \cite{D2}, to give a Quillen map
$\Spc'(k)_{\cT} \ra \Top$ we just need to give a functor $\Sm/k \ra
\Top$ which respects the above relations.  The functor we're
interested in is $X \mapsto X(\C)$, and this clearly preserves
relations (1) and (3).  In the case of the \'etale topology, the fact
that it preserves relation (2) is just Proposition~\ref{pr:genhyp};
the point is that if $p\colon E \map B$ is an \'etale cover, then
$p(\C)\colon E(\C) \map B(\C)$ 
satisfies the hypotheses of the inverse function theorem and hence
is locally split.

Since the \'etale topology is finer than the Nisnevich topology, there
is an obvious map $\Spc'(k)_{Nis} \ra \Spc'(k)_{et}$ (in essence, there
are more relations of type (2) for the \'etale topology).  So one also
gets a topological realization map $\Spc'(k)_{Nis} \ra \Top$ by
composition.
\end{proof}

It is possible to show that the functor $X \mapsto X(\C)$ takes
elementary distinguished squares \cite{MV} to homotopy pushout squares
of topological spaces.  Together with results of \cite{B}, this can
be used to give an alternative proof of the above theorem for the
Nisnevich topology.

\subsection{The Real case}
If we have a Real field $k\inc \R$, then the space $X(\C)$ comes
equipped with an action of the group $\Gal(\C/\R)=\Z_2$.  So we might
hope to compare $\Spc'(k)$ to a model category of $\Z_2$-equivariant spaces.

Recall that if $G$ is a finite group then there are two notions of
weak equivalence for $G$-spaces, called the \dfn{non-equivariant} and
\dfn{$G$-equivariant} equivalences.  
An equivariant map $X \ra Y$ is a non-equivariant equivalence if
it is a weak equivalence after forgetting the equivariant structure,
and it is a $G$-equivariant equivalence if 
$X^H \ra Y^H$ is a non-equivariant weak
equivalence for every subgroup $H\subseteq G$.
There are
associated $G$-equivariant and non-equivariant
model structures on the category of
$G$-spaces, which we will denote $\Gtop$ and $\Gtop_{\non}$.

If $p\colon E \ra B$ is an equivariant map which is also a covering
space (non-equivariantly), the map $\hocolim \ceck{E}_* \ra B$ is a
non-equivariant equivalence but not necessarily a $G$-equivariant
equivalence.  For instance, if $p$ is $G \ra *$ then the map $\hocolim
\ceck{E}_* \ra B$ is equal to $EG \ra *$.  So when we have a subfield
$k \inc \R$ the arguments given above show that the functor $X\mapsto
X(\C)$ induces a Quillen map $\Spc'(k)_{et} \ra \Zttop_{\non}$, but
not a Quillen map $\Spc'(k)_{et} \ra \Zttop$.  However, when we use
the Nisnevich topology something special happens.

\begin{lemma}
If $E \ra B$ is a Nisnevich cover of $k$-schemes, then
$E(\C)^{\Z_2} \ra B(\C)^{\Z_2}$ is locally split.
\end{lemma}

\noindent
For a counterexample to this in the case of \'etale covers, try
$\spec \C \ra \spec \R$.
  
\begin{proof}
First note that $X(\C)^{\Z_2}$ is homeomorphic to $X(\R)$ for any
scheme $X$ over $k$.  The map $p(\R): E(\R) \map B(\R)$ is surjective
by the defining property of Nisnevich covers; every $\R$-point in
$B$ lifts to $E$. 

By definition of \'etale covers, $p(\R)$ satisfies
the hypothesis of the inverse function theorem.  Since $p(\R)$ is 
surjective, it is locally split.
\end{proof}

\begin{thm}
There is a Quillen map $\Spc'(k)_{Nis} \ra \Zttop$ sending a smooth
$k$-scheme $X$ to $X(\C)$.  
\end{thm}

\begin{proof}
The argument exactly parallels the non-equivariant case in
Theorem \ref{th:Creal}, so the only nontrivial part is to show
that if $U_* \ra X$ is a Nisnevich hypercover then the map $\hocolim
U_*(\C) \ra X(\C)$ is a $\Z_2$-equivariant weak equivalence of
$\Z_2$-spaces.  The fact that it is a non-equivariant equivalence has
already been discussed in Theorem \ref{th:Creal}, because $U_* \ra
X$ is in particular an \'etale hypercover.  So we must consider what
happens when we take $\Z_2$-fixed points.

It is a fact that for any diagram $D$ of $G$-spaces ($G$ any finite
group) and any subgroup $H$ of $G$, one has $(\hocolim D)^H \he
(\hocolim D^H)$ (see Remark~\ref{re:fixedpts} below).  So we just need
to convince ourselves that $\hocolim \{U_*(\C)^{\Z_2}\} \ra
X(\C)^{\Z_2}$ is a non-equivariant weak equivalence.  But by the above
lemma one sees that $U_*(\C)^{\Z_2}$ is a generalized hypercover of
$X(\C)^{\Z_2}$, and so the result is an instance of
Proposition~\ref{pr:genhyp}.
\end{proof}

\begin{remark}
\label{re:fixedpts}
In the above proof we needed the fact that $(\hocolim D)^H$ is
weakly equivalent to $\hocolim (D^H)$.  This is well-known in
equivariant topology, but it's hard to find an actual reference.  We
give a brief sketch, for which we are grateful to Michael Mandell.

First of all, it clearly suffices to consider the case where all the
$D_i$ are cofibrant.  This means in particular that they are
Hausdorff.  We form $\hocolim D$ by first writing down the simplicial
replacement of the diagram, and then taking geometric realization.
Taking $H$-fixed points obviously commutes with the simplicial
replacement functor, so it suffices to worry about the geometric
realization part.  But one can check that if $X_*$ is a simplicial space
in which all $X_n$ are Hausdorff, then $\rea{X_*}^H$ is homeomorphic to
$\rea{X_*^H}$.  To do this, use the skeletal filtration on $\rea{X_*}$ and
the fact that $\rea{\Sk_n X_*}$ is obtained from $\rea{\Sk_{n-1} X_*}$ by
pushing out along a closed inclusion (this is one of the places where the
Hausdorff condition is needed).  Check that taking fixed-points
commutes with filtered colimits, and for Hausdorff spaces it also
commutes with pushouts along closed inclusions.
\end{remark}


\appendix

\section{Homotopy colimits for diagrams of non-cofibrant spaces}

Let $\Top$ denote the category of all topological spaces, with its
usual model category structure.  Given a diagram $D\colon I \ra \Top$,
the usual instructions for computing the homotopy colimit of $D$ are
(1) to apply a cofibrant-replacement functor to every object in the
diagram, and (2) to then use an explicit formula like that of
Bousfield-Kan \cite[Sec.~XII.2]{BK}.  This is the situation in an arbitrary
model category.  In this section we show that for the special case of
$\Top$, the first step of cofibrant-replacement is actually not
needed.  What we show is that no matter what formula one uses for
computing homotopy colimits---whether it is the Bousfield-Kan formula
or your favorite alternative---that formula always gives a homotopy
invariant construction in $\Top$, even without the
cofibrant-replacement step.  This fact seems not to be well known,
although it could be argued that the seeds lie there in the collective
subconscious of algebraic topologists.  In any case, for our purposes
here we need to bring it into the light of day.

The most useful way to formulate this result seems to be in model
category terms, as a comparison between the usual model structure on
$\Top$ and the Strom model structure, where everything is cofibrant.
See Theorem~\ref{th:hocolim}.

We would like to thank Phil Hirschhorn for helpful conversations about
the results in this section, in particular for his ideas on
removing an annoying $T_1$ separation condition.  The final form of
Lemmas~\ref{le:pushout-technical} and \ref{le:small} is something we
owe to him.

\bigskip

To begin with, we need the following

\begin{lemma}
\label{le:pushout}
Let $A \ra B$ and $X \ra Y$ be weak equivalences.
Given a diagram
\[ \xymatrix{
 A \times D^n \ar[d] & A\times S^{n-1} \ar[l]\ar[d]\ar[r] &X \ar[d] \\
 B \times D^n & B \times S^{n-1} \ar[l]\ar[r] & Y,}
\]
where the maps in the left-hand-square are the obvious ones, the
induced map from the pushout of the top row to the pushout of the
bottom row is also a weak equivalence.
\end{lemma}

Note that if $A$ and $B$ are cofibrant then this is an easy
consequence of left-properness for $\Top$, but we claim the
result in greater generality.

\begin{proof}
Let $X_A$ and $Y_B$ be the pushouts of the top and bottom rows
respectively, and write $f\colon X_A \ra Y_B$ for the map between
them.  We will produce a suitable cover of these spaces and use
Proposition \ref{pr:Gray}.

Let $U_B$ be the pushout of 
\[ \xymatrix{
 B\times (D^n-\{0\})  & B \times S^{n-1} \ar[r]\ar[l] & Y. }
\]
Write $D_\epsilon$ for $\{x \in D^n \colon\, \norm{x} < \epsilon\}$
(where $0 < \epsilon < 1$),
and let $V_B=B\times D_\epsilon$.
The spaces $U_B$ and $V_B$ clearly form an open cover of $Y_B$, and notice
that $U_B$ deformation-retracts down to $Y$.  The intersection $U_B\cap
V_B$ is equal to $B\times (D_\epsilon-\{0\})$.

The same definitions give us a cover $\{U_A,V_A\}$ of $X_A$, and it is
easy to check that $f^{-1}(U_B)=U_A$ and $f^{-1}(V_B)=V_A$.  So the
map $f^{-1}(V_B) \ra V_B$ is the map $A\times D_\epsilon \ra B \times
D_\epsilon$, which is a weak equivalence.  Similar reasoning shows
that $f^{-1}(U_B \cap V_B) \ra U_B \cap V_B$ is a weak equivalence.  Finally,
one argues that $f^{-1}(U_B) \ra U_B$ is a weak equivalence because it
deformation-retracts down to $X\ra Y$.  Proposition~\ref{pr:Gray} now shows
that $X_A \ra Y_B$ is a weak equivalence.
\end{proof}

We'll say that an inclusion $Y \inc Z$ is \mdfn{relatively $T_1$} if
given any open set $U$ in $Y$ and any point $z$ of $Z\backslash U$,
there is an open set $W$ of $Z$ such that $U\subseteq W$ and $z\notin
W$ (compare the similar definition from \cite[p.~50]{Ho}).  It follows
that if $E$ is any finite subset of $Z\backslash U$, one can find an
open set $W\subseteq Z$ which contains $U$ and doesn't intersect $E$.
Note that a space $X$ is $T_1$ precisely if all the inclusions $\{x\}
\inc X$ are relatively $T_1$.

\begin{lemma}
\label{le:pushout-technical}
Given a pushout diagram of the form
\[ \xymatrix{ A\times S^n \ar@{ >->}[d]\ar[r] & Y \ar@{ >.>}[d] \\
              A\times D^{n+1} \ar@{.>}[r] & Z,}
\]
the inclusion $Y \inc Z$ is relatively $T_1$.
\end{lemma}

\begin{proof}
Suppose given a point $z$ in $Z$ and an open $U$ in $Y$.  Either $z$ is
in $Y$ or else it is represented by a pair $(a,t)$ where $t$ is in the
interior of $D^{n+1}$.  The argument works the same for the two cases,
and so for convenience we'll assume the latter.

Pull back $U$ to $A\times S^n$ and express it as a union of rectangles
$V_i \times W_i$, where $V_i$ is open in $A$ and $W_i$ is open in
$S^n$.  Each $W_i$ can be fattened into an open subset $W'_i$ of
$D^{n+1}$ with the properties that $W'_i \cap S^n=W_i$ and 
$W'_i$ does not contain $t$.

Let $M$ be the union of the $V_i \times W'_i$; it is an open subset of
$A\times D^{n+1}$.  Let $N$ be the union of the images of $M$ and $U$
in $Z$.  One checks
that $N \cap Y=U$, and the pullback of $N$ to $A\times
D^{n+1}$ is $M$.  So $N$ is open in $Z$ and $N$ contains $U$, but $N$
does not contain $z$.
\end{proof}

The following lemma is well-known for closed inclusions of
$T_1$-spaces (see also \cite[Prop.~2.4.2]{Ho}).  The usual proof
still works in our case.

\begin{lemma}
\label{le:small}
Suppose that $Y_1 \inc Y_2 \inc \cdots$ is a sequence of relatively
$T_1$ inclusions and that $K$ is a compact space.  Then any map
$f\colon K\ra \colim Y$ factors through some $Y_k$.
\end{lemma}

\begin{proof}
Suppose the map does not factor through any $Y_k$.
By taking a subsequence of $Y$ if necessary, we can find a sequence
of points $k_1,k_2,\ldots$ in $K$ with the property that
$f(k_i)$ lies in $Y_i \backslash Y_{i-1}$.

Pick an $n$ and
set $V_n=Y_n$.  Next, choose an open set $V_{n+1}$ in $Y_{n+1}$ which
contains $V_n$ but doesn't contain $f(k_{n+1})$.  Then pick an open
set $V_{n+2}$ in $Y_{n+2}$ which contains $V_{n+1}$ but neither
$f(k_{n+1})$ nor $f(k_{n+2})$.  Continuing this process gives an
infinite sequence of opens, so their colimit $W_n$ is an open
subset of $\colim Y$.

As $n$ varies, the open subspaces $W_n$ form a cover of $\colim Y$.
But $f(K)$ is a compact subspace of $\colim Y$, and it is not covered
by any finite subcover.  This is a contradiction.
\end{proof}

We now need some machinery related to simplicial spaces.

\begin{defn}  
\label{de:freedeg}
A simplicial space $X_*$ is said to be \dfn{split}, or to have 
\dfn{free degeneracies}, if there exist subspaces $N_k\inc X_k$ such
that the canonical map
\[ \coprod_{\sigma} N_\sigma \ra X_k
\]
is an isomorphism.  Here the variable $\sigma$ ranges over all
surjective maps in $\Delta$ of the form $[k]\ra [n]$, $N_\sigma$
denotes a copy of $N_n$, and the map $N_\sigma \ra X_k$ is the one
induced by $\sigma^*\colon X_n \ra X_k$ (see \cite[Def.~8.1]{AM}).
\end{defn}

The idea is that the spaces $N_k$ represent the `non-degenerate' part
of $X_k$, sitting inside of $X_k$ as a direct
summand.  It is an easy exercise to check that if $X_*$ has free
degeneracies and all the $N_k$ are cofibrant spaces, then $X_*$ is
Reedy cofibrant in $s\Top$.  

If $X_*$ is any simplicial space, let $\Sk_n X_*$ be the simplicial
space equaling $X_*$ through dimension $n$ and equaling the degenerate
subspaces of $X_*$ in larger dimensions.  This is slightly different
than the $n$-truncated simplicial space $\sk_n X_*$.  There are maps
$\Sk_0 X_* \ra \Sk_1 X_* \ra \cdots$ and the colimit is $X_*$.  It
follows that $\rea{X_*}$ is equal to $\colim_n \rea{\Sk_n X_*}$, using
that geometric realization is a left adjoint (and this doesn't
require any assumptions on $X$, only hinging upon the fact that the
spaces $\del{n}$ are locally compact Hausdorff).  An important point is
that when $X_*$ has free degeneracies the space $\rea{\Sk_n X_*}$ is
obtained from $\rea{\Sk_{n-1} X_*}$ via the pushout diagram
\begin{equation}
\label{eq:sk}
 \xymatrix{   N_n \times \bd{n} \ar[d] \ar[r] 
            & \rea{\Sk_{n-1} X_*} \ar@{.>}[d]\\
                N_n \times \del{n} \ar@{.>}[r] & \rea{\Sk_n X_*}.}
\end{equation}

\begin{prop}
\label{pr:smallre}
Let $X_*$ be a simplicial space with free degeneracies.  If $K$ is
a compact space then any map $K \ra \rea{X_*}$ factors through some
$\rea{\Sk_n X_*}$.
\end{prop}

\begin{proof}
This is a direct application of Lemmas \ref{le:pushout-technical} and
\ref{le:small}, using the skeletal filtration of $\rea{X_*}$ and the
pushout square (\ref{eq:sk}).
\end{proof}

The following corollary is the crucial ingredient for
Theorem~\ref{th:hocolim}.  It is very similar to things in the
literature, notably \cite[Th.~11.13]{M1} and \cite[Lem.~A.5]{S2}.
May's result assumes the spaces are compactly-generated and Hausdorff,
and also that the realizations are simply-connected.  Segal's result
is more similar to ours, and the proofs follow the same pattern, but
he works with homotopy equivalences rather than weak equivalences.

\begin{cor}
\label{co:smallre}
If $X_* \map Y_*$ is a map of simplicial spaces with free degeneracies
such that $X_n \map Y_n$ is a weak equivalence for each $n$, then
$\norm{X_*} \map \norm{Y_*}$ is also a weak equivalence.
\end{cor}

\begin{proof}
For every $k$ and every basepoint $*$ of $X_0$, there is an
isomorphism
\[
\colim_n \pi_k (\norm{\Sk_n X_*}, *) \map \pi_k (\norm{X_*}, *)
\]
(and the same statement holds with $X_*$ replaced by $Y_*$).
This follows from Proposition \ref{pr:smallre}, taking $K$ to be a sphere.
Therefore, it suffices to show that $\rea{\Sk_n X_*} \map \rea{\Sk_n Y_*}$
is a weak equivalence.
Using induction, this follows from
the pushout square (\ref{eq:sk}) and 
Lemma~\ref{le:pushout}.
\end{proof}

Recall that the Strom model category is a model structure for
topological spaces, 
denoted $\Top^S$, in
which the weak equivalences are homotopy equivalences and the
cofibrations (resp., fibrations) are the Hurewicz cofibrations 
(resp., fibrations).
Note that all objects
are cofibrant in this structure.

\begin{prop}
The Strom model category is left proper and simplicial.
\end{prop}

\begin{proof}
Left properness is automatic when all objects are cofibrant
\cite[Cor.~11.1.3]{H}.  The simplicial action is of course given by
$A\tens K \iso A \times \norm{K}$.  To establish the simplicial structure
we use the
reductions outlined in \cite[Sec.~3]{D1}.  If $A\ra B$ is a
Hurewicz cofibration and $K\ra L$ is a cofibration of simplicial sets, then
$\norm{K} \ra \norm{L}$ is a closed cofibration and therefore the 
map
\[
A \tens L \coprod_{A \tens K} B \tens K \map B \tens L
\]
is a cofibration by \cite[Cor.~1]{L}.  If $A \trcof B$ is a
Hurewicz cofibration and a homotopy equivalence, then certainly $A\tens K
\ra B\tens K$ is still a homotopy equivalence.  And if $K\trcof L$ is
a trivial cofibration of simplicial sets then $\norm{K} \ra \norm{L}$
is actually a homotopy equivalence, hence $A\tens K \ra A\tens L$ is
also a homotopy equivalence.
\end{proof}

\begin{thm}
\label{th:hocolim}
Let $D\colon I \ra \Top$ be a diagram of spaces.  Then the homotopy
colimits of $D$ as computed in $\Top$ and $\Top^S$ have the same weak
homotopy type.
\end{thm}

\begin{proof}
In $\Top^S$, since all objects are cofibrant, we can compute $\hocolim
D$ by first taking the simplicial replacement of $D$ and then applying
the realization functor.
In $\Top$ we first apply a cofibrant-replacement functor to all
the objects in the diagram, and only then do we take simplicial
replacement and realize.  
Simplicial replacements always have free degeneracies (see \cite[Proof
of Lem.~2.7]{D2}), hence Corollary \ref{co:smallre} applies.
\end{proof}

\begin{remark}
Theorem~\ref{th:hocolim} also holds if one uses
the category of compactly-generated, weak Hausdorff spaces with its
usual model structure.  The same proofs work, with some extra caution
that the various colimits are what they're supposed to be.
\end{remark}

\bibliographystyle{amsalpha}

\end{document}